\documentclass[1pt]{article}
\usepackage{graphicx}

\usepackage{amsmath,amsthm}
\usepackage{amsfonts}

\newtheorem{theo}{Theorem}[section]
\newtheorem{prop}[theo]{Proposition}
\newtheorem{cor}[theo]{Corollary}

\newtheorem{defs}[theo]{Definitions}

\newcommand{\R}{\mathbb{R}}
\newcommand{\N}{\mathbb{ N}}
\newcommand{\C}{ \mathbb{C}}

\newcommand{\bc}{\begin{center}}
\newcommand{\ec}{\end{center}}
\newcommand{\be}{\begin{equation}}
\newcommand{\ee}{\end{equation}}
\newcommand{\bea}{\begin{eqnarray}}
\newcommand{\eea}{\end{eqnarray}}
\newcommand{\beq}{\begin{eqnarray*}}
\newcommand{\eeq}{\end{eqnarray*}}
\newcommand{\bv}{\left( \begin{array}{c} }
\newcommand{\ev}{\end{array} \right) }

\newcommand{\pei}{\rho_{ei}(T)}
\newcommand{\pea}{\rho_{e1}(T)}
\newcommand{\peb}{\rho_{e2}(T)}
\newcommand{\pec}{\rho_{e3}(T)}
\newcommand{\ped}{\rho_{e4}(T)}
\newcommand{\pee}{\rho_{e5}(T)}
\newcommand{\sei}{\sigma_{ei}(T)}
\newcommand{\sea}{\sigma_{e1}(T)}
\newcommand{\seb}{\sigma_{e2}(T)}
\newcommand{\sek}{\sigma_{e3}(T)}
\newcommand{\sed}{\sigma_{e4}(T)}
\newcommand{\see}{\sigma_{e5}(T)}

\newcommand{\Fp}{{\cal F_{+}}}
\newcommand{\Fm}{{\cal F_{-}}}

\newcommand{\KT}{{\cal K_{T}}}

\newcommand{\aaa}{\alpha}
\newcommand{\kkk}{\kappa}
\newcommand{\bbb}{\beta}
\newcommand{\ggg}{\gamma}
\newcommand{\lll}{\lambda}

\newcommand{\ooo}{\overline}

\title{Essential Spectra of Linear Relations}
\author{
     {\bf Diane Wilcox} \\
         School of Computational and Applied Mathematics\\
        University of the Witwatersrand,  Private Bag, Wits, 2050
        Johannesburg, {South Africa}	
}
\date{\today}

\begin{document}
\maketitle

\begin{abstract}
{Five essential spectra of linear relations are defined in terms of semi-Fredholm properties and the index.
Basic properties of these sets are established and the perturbation theory for semi-Fredholm relations is then applied to verify a
generalisation of Weyl's theorem for single-valued operators. We conclude with a M\"{o}bius transform spectral mapping theorem. }
\end{abstract}

\section{Introduction}

While the study of the spectrum of bounded linear operators
generalises the theory of eigenvalues of matrices, the essential
spectra of linear operators characterise the non-invertibility of
operators $\lambda - T$. The latter have been considered in terms
of two key related directions of investigation, namely the study
of the ascent and descent (as well as the nullity and defect) of
$\lambda - T$ and in terms of semi-Fredholm properties of $\lambda
- T$. Today there are several related definitions of essential
spectra and comprehensive reviews may be found in \cite{KoMu},
\cite{mm}, \cite{MbMu}, \cite{moll}, \cite{Rak} and \cite{Zel}. In
\cite{Lay1} , the refinements of the spectrum in terms of ascent
and descent were investigated in terms of states of operators,
using the terminology of  \cite{Tay2} (see also \cite{cross} for
the states of linear relations). On the other hand, the
perturbation theory of semi-Fredholm operators provides a more
general context for the early observations of H. Weyl, who showed
that limit points of the spectrum (i.e. all points of the
spectrum, except isolated eigenvalues of finite multiplicity) of a
bounded symmetric transformation on a Hilbert space are invariant
under perturbation by  compact symmetric operators (cf. Riesz and
Sz-Nagy \cite{RS}).

\vspace{0.3cm}
In this paper we apply the theory of Fredholm
relations to show that theory for essential
spectra of linear operators can be extended naturally to
linear relations. In particular, we extend preliminary results of Cross \cite{cross},
where  the set $\sigma_{e1}(\_)$ defined below is introduced.
The definitions in this paper are based on the classifications
given in Edmunds and Evans \cite{EdEv} for single-valued
operators.

\vspace{0.3cm}

 We commence with a recollection of some
preliminary properties required in the sequel.

\section{Semi-Fredholm Linear Relations}

We first clarify some notation and terminology. Let $X$ and $Y$ be
normed linear spaces, and let $B(X,Y)$ and $L(X,Y)$ denote the
classes of bounded and unbounded  linear operators, respectively,
from   $X$ into $Y$. A {\bf multivalued linear operator } $T:X
\rightarrow Y$ is a  set-valued map such that its graph $G(T) =
\{(x,y) \in X \times Y \ | \ y \in Tx \}$ is a linear subspace of
$X \times Y$. We  use the term{ \bf linear relation } or simply
relation, to refer to such a multivalued linear operator denoted
$T \in LR(X,Y)$ (cf. Arens \cite{arens} and Lee and Nashed
\cite{leena}). A relation  $T \in LR(X,Y)$ is said to be {\bf
closed} if its graph $G(T)$ is a closed subspace.
%We let $\tilde{X}$ denote the completion of the normed space $X$.
The {\bf closure}  of a linear relation $T$,
denoted  $\overline{T}$ is
defined in terms of its corresponding graph:
$\ G(\overline{T})  :=  \overline{G(T)} \subset X \times Y$.

The {\bf conjugate} $T'$ (cf \cite{cross}, III.1.1) of a linear relation $T \in LR(X,Y)$ is defined by
\bc $G(T') := G(-T^{-1})^{\bot} \subset Y' \times X'$ \ec
where  $[(y,x), (y', x')] := [x, x'] + [y, y'] = x'x + y'y. $
For $(y',x') \in G(T')$ we have $y'y = x'x$ whenever $x \in D(T).$

\newpage
Let $Q_T$, or simply $Q$, when there is no ambiguity about the relation $T$,
denote the natural quotient map
$\ Q_{_{\overline{T(0)}}}^Y:Y \rightarrow Y/\overline{T(0)}$
with kernel $\overline{T(0)}$.
For $x \in D(T)$ define $||Tx||$ by
\bc $||Tx|| := ||QTx||,$ \ec
and let the quantity $||T||$ be defined
\bc $||T|| := ||QT||.$ \ec
Clearly $QT$ is a single-valued linear operator.
It follows from the definition that $ ||Tx||  =  d(y, T(0)) $ for all $y \in Tx,$ and that
$||T|| = \sup\limits_{x \in B_{D(T)}} ||Tx||.$
The quantity $||T||$ is referred as the {\bf norm} of $T$, though we note that it is in fact a pseudonorm since
$||T|| = 0 $ does not imply $T = 0$.

\vspace{0.2cm}
A relation $T \in LR(X,Y)$  is said to be {\bf continuous} if for any neighbourhood
$V \subset R(T)$, the inverse image
$T^{-1}(V) := \{ u \in D(T) \ | \ V \cap Tu \neq \emptyset \}$ is a neighbourhood in $D(T)$, and $T$  is said to be {\bf open} if its inverse $T^{-1}$
is continuous.
It can be shown that
 $T$ is continuous if and only if $||T|| <  \infty$  (cf.  \cite{cross}, II.3.2).

\vspace{0.2cm}
The {\bf minimum modulus} of $T \in LR(X,Y)$ is the quantity
\bc $\ggg(T) := \sup\,\{ \lll \in \R \ : \ ||Tx|| \geq \lll d(x,N(T))$ for $x \in D(T) \},$
\ec
and {\em $T$ is open if and only if $\ggg(T) > 0$ }
(\cite{cross}, II.3.2). The quantity $\ggg(T)$ is related to the norm quantity by
$\ggg(T) = ||T^{-1}||^{-1}$.

%The injective component $\hat{T}$ of $T$ is defined by $\hat{T}:=Tq_{_{N(T)}}^{-1}\in LR(X/N(T), Y),$ where   %$q_{_{N(T)}}$ denotes the natural quotient from $X$ onto $X/N(T).$

\vspace{0.2cm}
The {\bf nullity} and {\bf deficiency } of a linear relation $T\in LR(X,Y)$ are defined respectively as follows:
\beq \aaa(T) &:=& \text{dim}N(T), \ \ \ \text{and} \\
        \bbb(T) &:=& \text{codim}R(T)\ := \ \text{dim}Y/R(T). \eeq

If either $\aaa(T) < \infty$ or  $\bbb(T) < \infty$, then the {\bf index} of $T$ is defined as follows:
\beq \kappa(T) &:=& \aaa(T)-\bbb(T), \eeq
where the value of the difference is computed as $\kappa(T):=\infty$ if $\aaa(T)$ is infinite and $\bbb(T)<\infty$ and $\kappa(T):=-\infty$ if $\bbb(T)$ is infinite and $\aaa(T)<\infty.$

\vspace{0.2cm}
If $X$ and $Y$ are Banach spaces and
$\ T:X \rightarrow Y\ $ is a closed single-valued operator, then
$T$ is said to be a {\bf Fredholm operator}, usually denoted $T \in \Phi(X,Y)$,
if  $R(T)$ is closed and both
$\aaa(T) < \infty$ and $\bbb(T) < \infty$;
$T$ is said to {\bf upper semi-Fredholm }, denoted $\ T \in \Phi_+(X,Y), \ $ if
$R(T)$ is closed and  $\aaa(T) < \infty;\ $ and $T$ is said to be
{\bf lower semi-Fredholm}, denoted $\ T \in \Phi_-(X,Y), \ $
if $R(T)$ is closed and $\bbb(T) <\infty.\ $

\vspace{0.1cm}

\begin{defs}
The {\bf essential resolvent sets}, $\ \pei$ for $i=1,2,3,4,5$,
of $T \in LR(X)$ are defined as follows:
\beq
\pea & := & \{\ \lambda \in \C\ |\ (\lambda - {T}) \in \Phi_{+} \cup \Phi_{-} \}  \\
\peb & := & \{\ \lambda \in \C\ |\ (\lambda - {T}) \in \Phi_{+} \} \\
\pec & := & \{\ \lambda \in \C\ |\ (\lambda - {T}) \in \Phi \} \\
\ped & := &\{\ \lambda \in \C\ |\ (\lambda - {T}) \in \Phi\
                                   and\ \kappa(\lambda - {T}) = 0 \} \\
\pee & := & \bigcup \ \rho_{e1}^{(n)}(T)\  where \
\rho_{e1}^{(n)}(T)\ \text{ is a component of } \; \rho_{e1} (T) \\
& &   \ \ \ \ \ \ \ \ \ \ \ \ \ \ \ \ \ \ and \ \
\rho_{e1}^{(n)}(T)\cap \rho(T)\ \neq\ \emptyset \eeq

The {\bf essential spectra}, $\ \sei,\ \ \
i=1,2,3,4,5, \ $  of $ \ T \in LR(X)\ $ are the respective
complements of the essential resolvents:

\vspace{0.3cm} \beq \sei := \ \C \setminus \pei,\ \  i=1,2,3,4,5.
\eeq

\newpage
We also define
\beq
\rho_{e2}'(T) & := & \{\ \lambda \in \C\ |\ (\lambda - {T}) \in \Phi_{-}\ \} \\
\sigma_{e2}'(T) & := & \ \C \setminus \rho_{e2}'(T) \eeq
\end{defs}

\vspace{0.5cm} Clearly we have that $\pei \supset \rho_{ej}(T)$
for $i<j<4$, and, thus, $\sei \subset \sigma_{ej}(T)$ for $i<j<4.$
We will see later that $\rho_{e4}(T) \supset \rho_{e5}(T)$.

\vspace{0.5cm}

For the rest of this section we recall a selection of results from Cross \cite{cross} which are used in the sequel .

\begin{prop}\label{cmpfr}
If $T\in LR(X,Y)$ is continuous with finite dimensional range, then $T$ is compact.
\end{prop}

\begin{prop}\label{sing2}
The following are equivalent:

(i) \ \ \  $T \not\in \Phi_+$.

(ii) \ \ There exists a non-precompact bounded subset $W$ of $D(T)$.

(iii) \ $T$ has a singular sequence.

\end{prop}

\begin{prop}\label{prev}
Let $T\in LR(X,Y)$ with $\gamma(T)>0$. Suppose $S \in LR(X,Y)$ satisfies $D(S) \supset D(T)$, $S(0) \subset \overline{T(0)}$
and $||S||<\gamma(T)$. Then $\aaa(T+S) \leq \aaa(T)$ and $\bar{\bbb}(T+S) \leq \bar{\bbb}(T)$.
\end{prop}

\vspace{0.1cm}
The next result is a general version of the so-called small perturbation theorem for linear relations.

\begin{prop}\label{ssper}
Let $S,T \in LR(X,Y)$. If $S(0) \subset \overline{T(0)}$ then
$\ \Delta(S) < \Gamma(T) \ \Rightarrow T+S \in \Phi_+, \ $ where \beq \Gamma(T) := \inf\limits_{M \in \cal{I}(D(T))}||T|_M||, \ \ \
 \Delta(S):=\sup\limits_{M \in \cal{I}(D(S))}\Gamma(S|_M), \eeq and $\cal{I}(X)$ denotes the collection of infinite dimensional subsets of $X$.
\end{prop}

\begin{prop}\label{ktwo}
Let $T \in \Phi(X,Y)$ and suppose $S \in LR(X,Y)$ satisfies $D(S) \supset D(T)$, $S(0) \subset \overline{T(0)}$
and $||S||<\gamma(T)$, then $\kappa(T+S) = \kappa(T)$.
\end{prop}

\begin{prop}\label{sper2}
Let $S,T \in LR(X,Y), D(S) \supset D(T)$ and let $T\in \Phi_-.$  \\
(a)  If $\dim R(S)< \infty$, then $T+S \in \Phi_-.$ \\
(b)  If $S$ is precompact, then $T+S \in \Phi_-.$\\
(c) If $||S||<\gamma(T')$, then $T+S \in \Phi_-.$
\end{prop}

\begin{prop}\label{ssp2}

(a) Suppose $T\in \Phi_+(X,Y)$ and $S \in LR(X,Y)$ is strictly singular. If $||S||<\infty$,  $D(S) \supset D(T)$, $S(0) \subset \overline{T(0)}$, then $\kappa(T+S) = \kappa(T)$.

(b) Suppose $T\in \Phi_-(X,Y)$ and $S \in LR(X,Y)$ is such that $S'$ is strictly singular. If $||S'||<\infty$,  $D(S) \supset D(T)$, $S(0) \subset \overline{T(0)}$, then $\kappa(T+S) = \kappa(T)$.

\end{prop}

\section{Properties of the Essential Spectra}

We begin this section by showing that the various essential
spectra are closed, and then illustrate some characteristic
properties. In the single-valued case, the set $ \
\bigcap\limits_{P \in \KT} \sigma{(T+K)} \ $ is  referred to as
the {\em Weyl essential spectrum}. Proposition ~\ref{charW} shows
that $ \ \sed \ $ can be characterised in terms  of the Weyl
essential  spectrum in the multivalued case as well (cf. Edmunds
and Evans \cite{EdEv}). We conclude this section by giving
properties of the quantities $\aaa (\lll - T),  \ $  $\bbb (\lll -
T)  \ $ and $\kkk (\lll - T)  \ $ for $\lll$ in the essential
spectra, and  deduce in  Proposition ~\ref{incl}   the inclusions
\beq \sea \ \subset \  \seb \ \subset \  \sek \ \subset \  \sed \
\subset \  \see \ \subset \ \sigma(T). \eeq Proposition
~\ref{per11} is included here for application in Proposition
~\ref{incl} and is based on the single-valued analogue given in
Goldberg \cite{gold}.

\vspace{0.2cm}
\begin{prop}\label{closed}
For $i=1,2,3,4,5, \  \sei$ is closed.
\end{prop}

{\sl PROOF}

Suppose $\ \lll \in \rho_{ei} (T), \ \ i=1,2,3,4, 5. \ $ Since $\
R(\lambda - T) \ $ is closed, it follows from the Open Mapping
Theorem (\cite{cross}, III.4.2), that $\gamma(\lambda - T) > 0$.
If $ \ \lll - T \in \Fp \ $ and  $\ |\mu| < \ggg (\lll-T), \  $  then
by Theorem ~\ref{ssper}, $\ \mu + \lll - T \in \Fp. \ $ Similarly,
if $ \ \lll - T \in \Fm \ $ and
 $\ |\mu| < \ggg (\lll-T'), \  $  then
by Theorem ~\ref{sper2}, $\ \mu + \lll - T \in \Fm. \ $ Thus,
$\pea, \ \peb$ and $\pec$ are open. Furthermore, by Theorem
~\ref{ktwo}, $\ \kappa ( \mu + \lll -T) = \kappa (\lll - T), \ $
i.e. $\ \ped \ $ is open. Since each component of $\ \pea\ $ is
open, so is $\ \pee. \ $

\begin{prop}
Let $T \in LR(X)$ . Then

\vspace{0.2cm} (a) $\sigma_{ei}(T') = \sei$ for $i=1,3,4,5$

\vspace{0.1cm} (b) $\sigma_{e2}(T') = \sigma_{e2}'(T)$
\end{prop}

{\sl PROOF}

(a) Suppose $\ \lll \in \rho_{ei} (T), \ \ i=1,3,4. \ $ By
\cite{cross}, III.7.2,  $\ \alpha(\lambda - T') = \beta(\lambda
- T)\  $
 since  $ \ R(\lambda - T) \ $  is closed.
By the Closed Range Theorem (\cite{cross}, III.4.4), $ \ R(\lambda - T') \ $ if
and only if  $ \ R(\lll - T) \ $ is closed and, since $\lll-T$ is
open,
 $\ \bbb(\lambda - T') =
\aaa(\lambda - T). \ $ Thus,  the result holds for $i=1,3$ and
$4$. Since  $\ \pea = \rho_{e1}(T') \ $ and $ \ \rho(T) =
\rho(T'), \ $ it follows that $ \ \rho_{e1}^{(n)}(T') =
\rho_{e1}^{(n)}(T), \ $ i.e. the result holds for $i=5$.

\vspace{0.3cm} (b) follows from the reasons given in (a).

\begin{prop}
$\lambda \in \seb $ if and only if $ \lambda - T$ has a singular
sequence.
\end{prop}

{\sl PROOF}

Since $\lambda \in \seb$ if and only if $ \lambda - T \notin \Fp$,
the result follows from Theorem ~\ref{sing2}.

\begin{prop}\label{charW}

\beq \sed\ =\ \bigcap\limits_{K \in \KT} \sigma{(T+K)}, \eeq

where $ \  \KT := \{K \in LR(X)  \ | \ K \ \ is  \ \ compact  \ \
and \ \ K(0) \subset \ooo{T(0)} \ \ \}. \ $

\end{prop}

{\sl PROOF}

We show first that $\ \sed\ \subset\ \bigcap_{K \in \KT}\
\sigma{(T+K)}. \ $ Suppose $ \ \lambda \notin \bigcap_{K \in \KT}\
\sigma{(T+K)}. \ $ Then there exists $ \  K \in \KT \ $ such that
$ \ \lambda \in \rho(T+K). \ $ Thus $\ \lambda \in \rho_{e4}(T+K).
\ $ By Propositions ~\ref{ssper} and ~\ref{sper2},  $  \lll - T =
\lll - T - K + K  \in \Phi, \ $ and by Theorem ~\ref{ssp2},
 \beq  \kkk(\lll - T) = \kkk(\lll - T - K + K) = \kkk(\lll-T-K).  \eeq
Thus, $\ \lambda \in \rho_{e4}(T), \ $ i.e.  $\ \lambda \notin
\sed$.

\vspace{0.3cm} Conversely, suppose $ \ \lambda \in \rho_{e4}(T) .
\ $ Then $\ R(\lambda - T) \ $ is closed, and $ \ \alpha(\lambda -
T) = \beta(\lambda - T) = n, \ $ say. Let $\ \{x_{1}, \ldots ,
x_{n} \}  \ $ and $ \  \{y_{1}', \ldots , y_{n}' \} \ $ be bases
for $ \ N(\lambda - T) \ $ and $\ R(\lambda - T)^{\bot} =
N(\lambda - T'), \ $ respectively. Choose $x_{j}' \in X'$ and
$y_{j} \in X,\ j=1, \ldots , n \ $ such that
\beq x_{j}'x_{k} & = & \delta_{jk},  \ \ and \\
y_{j}'y_{k} & = & \delta_{jk}, \eeq where $\delta _{jk} = 0$ if $j
\neq k$ and $\delta _{jk} = 1$  if $j=k, \ $ and define $K \in
LR(X) $ as follows: \beq  Kx \ := \
\sum\limits_{k=1}^{n}\,(x_{k}'x)\,y_{k}, \ \ \ \ x \in X \eeq Then
$ \ \dim\,R(K) < \infty\  $and \beq ||Kx|| \ \leq \
(\sum\limits_{k=1}^{n}||x_{k}'||\,||y_{k}||)\,||x||. \eeq By
Proposition ~\ref{cmpfr}, it follows that  $K$ is a compact
operator. By Propositions ~\ref{ssper} and ~\ref{sper2},  it
follows that  $\ \lll - (T+K) \in \Phi \ $ and by Theorem
~\ref{ssp2}, $ \ \kappa (\lll- (T+K)) = \kappa (\lll - T). \ $

\vspace{0.5cm} Without loss of generality, assume $\lambda = 0$.
Now if $x \in N(T),$ then $x = \sum\limits_{k=1}^{n}a_{k}x_{k}$
and $x_{j}'(x) = a_{j},\ 1 \leq j \leq n$. On the other hand, if
$x \in N(K)$, then $x_{j}'(x) = 0$. Thus $N(T) \cap N(K) = 0.$

\vspace{0.5cm} Similarly, if $y \in R(K),$ then $y =
\sum\limits_{k=1}^{n}a_{k}y_{k}$ and $y_{j}'(y) = a_{j},\ 1 \leq j
\leq n, \ $ and if $y \in R(T)$, then $y_{j}'(y) = 0$. Thus $R(K)
\cap R(T) = 0$

\vspace{0.5cm} Next, suppose $x \in N(T+K)$. Then $Tx = -Kx +
T(0). \ $ It follows from the argument above, that $Tx = T(0),$
i.e. $x \in N(T)$. Thus, $\ x = \sum\limits_{k=1}^{n}a_{k}x_{k}$
and $\ x_{k}'(x) = a_{k},\ 1 \leq k \leq n$. Since $\ Kx =
\sum\limits_{k=1}^{n}(x_{k}'x)y_{k} = 0$, it follows that
$x_{k}'(x) = 0,\ 1 \leq k \leq n$, and hence $x=0. $ Thus,
$\alpha(T+K) = 0 = \beta(T+K), \ $ i.e. $0 \in \rho_{e4}(T+K)$ .

\begin{prop}\label{per11}
Suppose $ \ T\in \Phi_+ \cup
\Phi_- \ $  and $S\in LR(X,Y) $  satisfies $\ D(S) \supset D(T), \
$ $S(0) = \ooo{S(0)} \subset \ooo{T(0)}, \ $ and $ \ ||S||<
\ggg(T). \ $ Then $\exists \ \nu > 0$ such that $\ \aaa (T+ \lll
S)$ and $\bbb (T+ \lll S)$ are constant in the annulus $0 < |\lll
| < \nu$.
\end{prop}

{\sl PROOF}

We first assume $\aaa (T) < \infty$. Let $\lll \neq 0$ and let $x
\in N(T + \lll S).$ Then \beq Tx \supset - \lll Sx, \eeq whence
\beq
Sx & \subset &  R(T) =:R_1, \ \  and \\
 x & \in & S^{-1}R_1 = : D_1. \eeq
Thus \beq
-\lll Sx & \subset &  Tx \subset TD_1 = : R_2,  \ \ and \\
 x & \in & S^{-1}R_2 = : D_2. \eeq
Proceeding in this way, we obtain \beq R_{k+1} : = TD_k, \ \ where
\ \  D_k := S^{-1}R_k.  \eeq Clearly \bc $R_1 \supset  R_2 \supset
\ldots \ \ $
     and $\ D_1 \supset  D_2 \supset \ldots \ $ \ec
It follows from the construction of these sequences of subspaces
that \bea N( T + \lll S) \ \subset  \
\bigcap\limits_{k=1}^{\infty}\,D_k. \label{q6} \eea By induction,
we have that $R_n$ are closed subspaces of $Y,$ and  $D_n$ are
relatively   closed subspaces of $D(S)$: from the hypothesis,
$R_1$ is closed, and, hence, since $S$ is continuous, and $S(0)$
is closed, $D_1$ is relatively  closed in $D(S)$; if  $R_k$ and
$D_k$ are closed and relatively closed, respectively, then, since
$ \ T|_{D_k} \in \Phi_{+} \cup \Phi_-, \ $ it follows that
$R_{k+1} = TD_k$ is closed, and, since $S$ is continuous, and
$S(0)$ is closed,
 $D_{k+1} = S^{-1}R_{k+1}$ is relatively closed in $D(S)$.

\vspace{0.3cm} Define
\beq X_1 \ & := & \  \bigcap\limits_{k=1}^{\infty}\,D_k, \ \  and \\
     Y_1 \ & := & \  \bigcap\limits_{k=1}^{\infty}\,R_k.
\eeq Then, by the definitions of $R_k$ and $D_k$, it follows that
    \bc $TX_1 \subset Y_1 \ $ and $\ SX_1 \subset Y_1.$\ec
Now define $T_1$ and $S_1$ by : \beq T_1 := T|_{D(T) \cap X_1}, \
\ \ and \ \  \
    S_1 := S|_{D(T) \cap X_1}.
\eeq Then $R(T_1) \subset Y_1$ and $R(S_1) \subset Y_1, $ and
since $T$ is closed  and $X_1$ is relatively closed in $D(S)$ and
hence also in $D(T)$,
 $T_1$ is a closed relation.
To see  that $T_1$ is surjective, let $y \in Y_1 =
\bigcap\limits_{n=1}^{\infty}TD_n$. Then for each $n \geq 1$,
there exists $x_n \in D_n$ such that $y \in Tx_n$. Since  $\aaa
(T) < \infty$ and $D_n \supset D_{n+1}$, there exists  $ k_0$ such
that for $k \geq k_0,$ \beq N(T) \cap D_{k_0} = N(T) \cap D_k,
\eeq and for $\ x_k \in D_k,  \ $ and $ \ x_{k_0} \in D_{k_0}, $
\beq
 x_k - x_{k_0} \in N(T) \cap D_{k_0} = N(T) \cap D_k
                 \subset D_k.
\eeq From this it follows that \beq x_{k_0}  \in \bigcap\limits_{k
\geq k_0}\,D_k = X_1,  \ \ and \ \
 y  \in  Tx_{k_0}.
\eeq i.e.  $T_1$ is surjective. By the Open Mapping Theorem
(\cite{cross}, III.4.2), $T_1$ is open.

\vspace{0.3cm} By  Theorem ~\ref{prev}, Propositions ~\ref{ssper}
and ~\ref{sper2}, and by Theorem ~\ref{ktwo} , $\exists \ \nu > 0$
such that for $ |\lll |  < \nu$  we have

\vspace{-0.5cm}

\bea
 \kappa (T+ \lll S) & = & \kappa (T)  \label{q3}.  \eea

\vspace{-0.5cm}

Since

\vspace{-0.5cm}

\bea  \bbb (T_1+ \lll S_1) \ \leq \ \bbb (T_1)  \ = \
\bar{\bbb}(T_1) \ = \ 0 \label{q4}, \eea it follows that $\  \bbb
(T_1+ \lll S_1) = 0, \ $ and hence \bea
  \aaa (T_1+ \lll S_1)\  = \ \kappa (T_1+ \lll S_1) & = & \kappa (T_1) \  = \ \aaa (T_1) \label{q5}.
\eea By ( ~\ref{q6}), it follows that for $\lll \neq 0, $ \be N(T+
\lll S) = N(T_1+ \lll S_1). \label{q7} \ee In particular, $ \ \aaa
(T+ \lll S) = \aaa (T_1+ \lll S_1)$. By (~\ref{q3}), ( ~\ref{q4}),
(~\ref{q5}) and  (~\ref{q7}) it follows that $\aaa (T+ \lll S)$
and $\bbb (T+ \lll S)$ are constant in the annulus $0 < |\lll |  <
\nu$.

\vspace{0.3cm} If $\aaa(T) = \infty,$ then $\bbb(T) < \infty$, and
the result is obtained by passing to the conjugates.

%\begin{prop}
%Let $X$ and $Y$ be complete, and suppose
%$ \ T + \lll S \Phi_+ \cup \Phi_- \ $  and $S\in LR(X) $ is continuous and satisfies %$\ D(S) \supset D(T)  \ $ and
%$S(0) = \ooo{S(0)} \subset \ooo{T(0)}. \ $
%Then
%then
%$\aaa (T+ \lll S)$ and $\bbb (T+ \lll S)$ have constant values, $n_1$ and $n_2,$
%respectively, except perhaps at isolated points where
%\bc $\aaa (T+ \lll S) > n_1$  and $\bbb (T+ \lll S) > n_2.$
%\ec
%\end{prop}
%
%\begin{prop}\label{const}
%Let $X$ and $Y$ be complete and let $S$ be continuous with $D(T) \subset D(S)$.
%If $\rho_{e1}^{(n)}(T)$ is a component of $\rho_{e1}(T),$ then
%$\aaa (T+ \lll S)$ and $\bbb (T+ \lll S)$ have constant values, $n_1$ and $n_2,$
%respectively, except perhaps at isolated points where
%\bc $\aaa (T+ \lll S) > n_1$  and $\bbb (T+ \lll S) > n_2.$
%\ec
%\end{prop}

\begin{prop}\label{const}
%Let $X$ be complete and let $T \in LR(X)$ be closed.
If $\ \rho_{ei}^{(n)}(T)\ $ is a component of $\ \rho_{ei}(T),\ i=1,2,3,
\ $ then $\ \aaa (\lll - T) \ $ and $\ \bbb (\lll - T)\ $ have
constant values, $n_1$ and $n_2,$ respectively, $n_1, n_2 \in \N
\cup \{\infty\}$, except perhaps at isolated points where \bc
$\aaa (\lll - T) > n_1$  and $\bbb (\lll - T) > n_2.$ \ec
\end{prop}

{\sl PROOF}

We first prove the result for the quantities $\aaa(\lll - T)$.
Since any component of an open set in $\C$  is open, we have that
$\rho_{ei}^{(n)}(T) \  $  are  open sets. We first consider the
case $\rho_{e1}^{(n)}(T)$. If $\aaa(\lll-T) = \infty$ for all $
\lll \in \rho_{e1}^{(n)}(T)$, then we are done. Now suppose
$\aaa(\lll-T) < \infty$ for some  $ \lll \in \rho_{e1}^{(n)}(T)$,
define $\aaa( \lll) := \aaa (\lll - T)$, and choose $\lll_0$ such
that  $\aaa(\lll_0) = n_1$ is the smallest non-negative integer
attained by  $\aaa( \lll)$ on $\rho_{e1}^{(n)}(T)$. Suppose $\aaa(
\lll') \neq n_1$ for some $ \lll ' . $ Since $\rho_{e1}^{(n)}(T)$
is connected, there exists an arc $\Lambda$ in
$\rho_{e1}^{(n)}(T)$ with endpoints $\lll_0$ and $\lll'$. Since $
\lll - T \in \Phi_+ \cup \Phi_- \ $, it follows from
 Proposition ~\ref{per11}  that for each
$\mu \in \Lambda$ there exists an open ball $B_{\mu}$ contained in
$\rho_{e1}^{(n)}(T)$ such that   $\aaa( \lll)$ is constant on
$B_{\mu} \setminus \{\mu\}$. Since $\Lambda$ is compact, there
exists a finite set of points $\lll_1, \lll_2, \ldots, \lll_n =
\lll'$ such that $\ B_{\lll_0}, B_{\lll_1}, \ldots, B_{\lll_n}$
cover $\Lambda$, and, for $0 \leq i \leq n-1,$ \be B_{\lll_i} \cap
B_{\lll_{i+1}} \neq \emptyset.   \label{q8} \ee It follows from
Theorem ~\ref{prev} that $\aaa(\lll) \leq \aaa(\lll_0)$ for $\lll$
sufficiently close to $\lll_0$. Thus, since $ \aaa(\lll_0)$ is the
minimum value attained by  $\aaa(\lll)$ on $\rho_{e1}^{(n)}(T)$,
it follows that $\aaa(\lll) = \aaa(\lll_0)$ for $\lll$
sufficiently close to $\lll_0$. Since $\aaa(\lll)$ is constant for
all $\lll \neq \lll_0$ in $B_{\lll_0}$, this constant must be
$\aaa(\lll_0)$. Similarly $\aaa(\lll)$ is constant on $B_{\lll_i}
\setminus \{\lll_i\}$ for $1 \leq i \leq n. $ Thus, by (~\ref{q8})
that $\aaa(\lll) = \aaa(\lll_0)$ for all $\lll \in B_{\lll'}
\setminus \{\lll'\}$ and $\aaa(\lll') > n_1$.

\vspace{0.3cm} To see that the result holds for $\bbb (\lll - T)$,
we pass to the conjugate of $T$  and apply the above, and the
equality \bc $\aaa (\lll -  T') = \bbb (\lll - T).$ \ec

The proofs for $\rho_{e2}^{(n)}(T)$ and $\rho_{e3}^{(n)}(T)$ are
similar.

\begin{prop}\label{pee}
$\lll \in \pee$ if and only if $\lll \in \ped$ and a deleted
neighbourhood of $\lll$ lies in $\rho(T)$.
\end{prop}

{\sl PROOF}

Suppose $\lll \in \pee$. Then, by definition,  $\lll$ lies in a
component $\rho_{e1}^{(n)}(T)$ of $\rho_{e1}(T)$ which intersects
$\rho(T)$. Let $C$ be such a component. Clearly $C \cap \rho(T)$
is open.

\vspace{0.3cm} Since $\mu \in C \cap \rho(T)$ implies $\aaa(\mu -
T) = \bbb(\mu - T) = \kappa(\mu - T)  = 0$, it follows  from
Theorem ~\ref{ktwo} that $\kappa(\lll - T)  = 0$ for $\lll \in C$
when $\lll$ is sufficiently close to $\mu$, and, hence for all
$\lll \in C$ . Applying Proposition ~\ref{const}, we see that
$\aaa(\lll - T) = \bbb(\lll - T) = 0$ for all except some isolated
points, say $\lll_j$ where $\aaa(\lll_j - T) > 0$ and $\bbb(\lll_j
- T) > 0$. Thus if $\lll \in \pee$, then either $\lll \in \rho(T)$
or $\lll$ is one of these isolated points in $\rho_{e4}(T)$.

\vspace{0.2cm} Clearly the converse is true.

\begin{cor}\label{pee2}
If $\rho_{e4}(T)$ is connected and $\rho(T) \neq \emptyset$, then
$\rho_{e5}(T) = \rho_{e4}(T).$
\end{cor}

{\sl PROOF}

Since $\rho(T) \subset  \rho_{e4}(T)$, it follows from the
hypothesis and Proposition ~\ref{const}   that $\aaa(\lll - T) =
\bbb(\lll - T) = 0$ for all $\lll \in \rho_{e4}(T)$ except perhaps
at isolated points, \ i.e.  a deleted neighbourhood of $\lll$ lies
in $\rho(T)$. The result follows from Proposition ~\ref{pee}.

\begin{prop}\label{incl}
\beq \sea \ \subset \  \seb \ \subset \  \sek \ \subset \  \sed \
\subset \  \see \ \subset \ \sigma(T) \eeq
\end{prop}

{\sl PROOF}

Clearly \beq \pea \ \supset \  \peb \ \supset \  \pec \ \supset \
\ped. \eeq The remaining inclusions follow from Proposition
~\ref{pee}.

\begin{prop}
The index is constant in each connected component
$\rho_{ek}^{(n)}(T) $ of
 $\rho_{ek}(T), k=1,2,3,4,5.$
\end{prop}

{\sl PROOF}

Clearly the result holds for $\rho_{e4}^{(n)}(T)$, and it follows
from Proposition ~\ref{pee} that the result hold for
$\rho_{e5}^{(n)}(T)$.

\vspace{0.3cm} Let $\lll$ and $\lll'$ be distinct points in
$\rho_{ek}^{(n)}(T), k=1,2,3$. Let $\Lambda$ be an arc in
$\rho_{ek}^{(n)}(T)$ with endpoints $\lll$ and $\lll'$. By Theorem
~\ref{ktwo},  there exists $\epsilon > 0$ such that $\kappa(\mu -
T) = \kappa(\lll -T)$ for any $\mu$ such that $|\mu - \lll| <
\epsilon$. Clearly the open balls $B(\lll), \  \lll \in \Lambda\ $
cover $\Lambda$. Since $\Lambda$ is compact, a finite  number of
these balls suffices to cover $\Lambda$. Since each of these balls
overlap, it follows that $ \kappa (\lll - T) = \kappa (\lll' -
T)$.

\section{Perturbation of the Essential Spectra}

We now apply perturbation theorems for semi-Fredholm
relations to verify the stability properties of the essential
spectra under small and compact perturbation. In particular we
arrive at  generalisations of Weyl's theorem for  linear
operators to a relatively compact case \cite{BH} . First we recall Propositions ~\ref{gtg} to ~\ref{rbcl} which are proved in Cross \cite{cross}.

\begin{prop}\label{gtg}
Let $T \in LR(X,Y)$ and let $G=G_T$ denote the graph operator of $T$,
i.e. $G_T$ is the identity injection of $X_T$ into $X$ ($G_Tx=x$) and $X_T$ is the vector space $D(T)$ endowed
with the norm $||x||_T:= ||x|| +||Tx||$ for $x \in D(T)$. Then $TG$ is open if and only if $T$ is open
and \beq \gamma(TG) = \frac{\gamma(T)}{1+\gamma(T)}, \ \text{provided}\ \ T \neq 0,\eeq
with the cases $\frac{\infty}{\infty}:=1$ and $\gamma(TG) :=\infty$ if $T=0$.
\end{prop}

\begin{prop}\label{rbcl2}
 The norms $||\_||_T$ and $||\_||_{\lambda-T}$ are equivalent.
\end{prop}

\begin{prop}\label{rbcl}
Let $T\in LR(X,Y)$ and suppose $S\in LR(X,Y)$ satisfies $\ D(S) \supset \overline{D(T)} \ $ and
$S(0)  \subset {T(0)}, \ $ and is $T$-bounded with $a,b>0, \ b<1$   such that for $x \in D(T)$, $||Sx||\leq a||x||+b||Tx||.$

(a)  The norms $||\_||_T$ and $||\_||_{T+S}$ are equivalent.

(b) If $X$ and $Y$ are complete and $T$ is closed, then $T+S$ is closed.
\end{prop}

\begin{theo}\label{Weyl}
Let $T \in LR(X)$ be closed and suppose $S \in LR(X)$ is $T-compact$ with $T-$bound $b < 1$ \cite{BH}, and
$D(S) \supset \ooo{D(T)}$ and $S(0) \subset T(0)$. Then for
$i=1,2,3,4$ \bc $\sigma_{ei}(T+S) = \sigma_{ei}(T).$ \ec If
additionally $\rho_{e4}$ is connected and neither $\rho(T)$ nor
$\rho(T+S)$ are empty, then \bc $\sigma_{e5}(T+S) =
\sigma_{e5}(T).$ \ec
\end{theo}

{\sl PROOF}

By Corollary ~\ref{rbcl2},  the norms $||\_||_{_T}$ and
$||\_||_{_{\lll - T}}$ are equivalent and hence,  $S$ is $(\lll -
T)-compact$. Let $\ G_{_{\lll-T}}$ denote the graph operator from
space $\ X_{_{\lll - T}} : = \ (X,  ||x||_{_{\lll-T}}) \ $ into
$X$. Suppose $\lll-T \in \Phi_{\pm}$. Clearly  $R(TG_{_{\lll-T}})
= R(T), $ and as subsets of the set $X$, we have
$N(TG_{_{\lll-T}}) = N(T)$. By Proposition ~\ref{gtg},
$(\lll-T)G_{_{\lll-T}}$ is open, and hence $(\lll-T)G_{_{\lll-T}}
\in \Phi_{\pm}$. Thus, by Propositions ~\ref{ssper} and
~\ref{sper2},  it follows that $(\lll-T) - S = \lll - (T+S) \in
\Phi_{\pm}$ and
 by Theorem ~\ref{ssp2}, $\kappa(\lll - (T+S)) = \kappa(\lll - T)$.

\vspace{0.3cm} On the other hand, suppose $\lll-(T+S) \in
\Phi_{\pm}$. By the equivalence of the norms $||\_||_T$ and
$||\_||_{\lll - (T+S)}$ (Proposition ~\ref{rbcl} and Corollary
~\ref{rbcl2}), it follows that $S$ is $(\lll - (T+S))-compact$.
Arguing as before, it follows that $\lll-T \in \Phi_{\pm}$ and $
\kappa(\lll - T) = \kappa(\lll - (T+S))$.

\vspace{0.3cm} Thus, $\rho_{ei}(T+S) = \rho_{ei}(T)$ for
$i=1,2,3,4$. It follows from the additional hypotheses, Corollary
~\ref{pee2}, and what has just been proved that \bc $\rho_{e5}(T)
= \rho_{e4}(T)  = \rho_{e4}(T+S) = \rho_{e5}(T+S). $  \ec

\section{Functions of the Essential Spectra}

The  M\"{o}bius transform, $\eta(\lambda) = (\mu - \lambda)^{-1}$, is a topological homeomorphism from $\ \C \cup \{\infty\}$, endowed with the usual topology, onto itself. Theorem ~\ref{Moeb2} below is analogous to the Theorem on the  M\"{o}bius transform of the spectrum in Cross \cite{cross}. For its proof, we first recall the following index theorem:

\begin{prop}\label{ind1}
Let $T\in LR(X,Y)$ and $S \in LR(YZ)$. Suppose $D(S) =Y$ and that $T$ and $S$ have finite indices. Then
\beq \kappa(ST) = \kappa(T) + \kappa(S) - \dim(T(0) \cap N(S) ). \eeq
\end{prop}

\begin{theo}\label{Moeb2}
Let $T \in LR(X)$ be closed. Suppose $\mu
\in \rho(T)$. Then for $i=1,2,3,4,5$ \bc $\lll \in  \sigma_{ei}(T)
\Leftrightarrow
    (\mu - \lll)^{-1} \in \sigma_{ei}(T_{\mu}).$
\ec
\end{theo}

{\sl PROOF}

Let  $S := (\mu - \lll)((\mu - \lll)^{-1} - T_{\mu})$.  It can be shown  that  $\lll -
T = S(\mu - T)$ (\cite{cross}, IV.4.2).   Since T is closed, so is $\lll - T$, and
since $R(\mu - T) = X$ it follows that \be    R(\lll - T) = R(S).
\label{q9} \ee Since $T_{\mu}$ is single valued, \bc $\aaa(\lll -
T)  = \dim T_{\mu}S^{-1}(0) \leq \dim S^{-1}(0) = \aaa(S).$ \ec
Thus, $S \in \Phi_{\pm}$ implies  that $\lll - T \in \Phi_{\pm},\
$ i.e.
 $(\mu - \lll)^{-1} \in \rho_{ei}(T_{\mu})$ implies that
$ \lll \in  \rho_{ei}(T)$ for $i=1,2,3. \ $ Applying Proposition
elementary algebra for linear relations (\cite{cross}, I.4.2) we have
\beq (\mu-T)S \ & =& \ (\mu-T)(\mu - \lll)((\mu - \lll)^{-1} - T_{\mu}) \\
& = & \ (\mu-T) - (\mu-\lll)(\mu-T)(\mu-T)^{-1} \\
& = & \ (\mu-T) - (\mu-\lll)(I + (\mu-T)(\mu-T)^{-1} - (\mu-T)(\mu-T)^{-1} ) \\
& = & \ \lll - T + (\mu - \lll)(TT^{-1} - TT^{-1}) \\
& = & \ \lll - T. \eeq Thus, since $ \kkk(\mu-T) $and $\kkk(S)$
are finite and $D(S) = X$, it follows from Proposition ~\ref{ind1}
that \be  \kappa(\lll - T) = \kappa(S) + \kappa (\mu - T) -
\dim(S(0) \cap N(\mu - T)). \label{q10} \ee In particular, if
$(\mu - \lll)^{-1} \in \rho_{e4}(T_{\mu})$ then $\kappa(S) = 0$,
and, since $\mu \in \rho(T)$,\ we have \linebreak $\kappa(\mu - T)
= 0  = \aaa(\mu - T)$. Thus $\kappa(\lll - T) = 0, \ $ i.e. $\lll
\in \rho_{e4}(T)$. Applying Proposition ~\ref{pee}, it follows
that the forward implication also holds for $i=5$.

\vspace{0.3cm} For the reverse implication, it follows from
(~\ref{q9}) that if $\lll - T \in \Phi_{-}$, then $S \in \Phi_{-},
\ $ i.e.  \linebreak $(\mu - \lll)^{-1}-T_{\mu} \in \Phi_{-}$. Now
suppose $\lll - T \in \Phi_{+}$. Then there exists a finite codimensional subset$ M$ of $D(\lll - T)$ such that $(\lll - T )|_M$ is injective. As in
\cite{cross} IV.4.2,  it follows that $S|_M$ is injective, and
hence $\aaa(S) < \infty$. Thus, $S \in \Phi_{+}$, and consequently
$(\mu - \lll)^{-1} - T_{\mu} \in \Phi_{+}$. We have \bc $\lll \in
\rho_{ei}(T) \ \Rightarrow \ (\mu - \lll)^{-1} \in
\rho_{ei}(T_{\mu})\ $ for $i=1,2,3.$ \ec Now if $\lll \in
\rho_{e4}(T)$ then $\kappa (\lll - T) = 0$, and since $\aaa (\mu -
T) = \kappa (\mu - T) = 0$ it follows from (~\ref{q10}) that  $0 =
\kappa(S) = \kappa ((\mu - \lll)^{-1} - T_{\mu})$. Thus $(\mu -
\lll)^{-1} \in \rho_{e4}(T_{\mu})$. Another application of
Proposition ~\ref{pee} shows that the converse is true for $i=5$.

\begin{theo}
Let $X$ be complete and let $T, S \in LR(X)$ be closed.

Suppose $\mu \in \rho(T) \cap \rho(S)$ and $T_{\mu} - S_{\mu}$ is
compact. Then for $i=1,2,3,4$ \bc $\sigma_{ei}(S) =
\sigma_{ei}(T).$ \ec If additionally $\rho_{e4}(S)$ is connected
then equality holds for $i=5$ as well.
\end{theo}

{\sl PROOF}

For $i=1,2,3,4$ it follows from Theorem  ~\ref{Moeb2},  that \beq
\lll \in \sigma_{ei}(T) & \Leftrightarrow & (\lll - \mu)^{-1} \in
\sigma_{ei}(T_{\mu}),\eeq and \beq \lll \in \sigma_{ei}(S) &
\Leftrightarrow & (\lll - \mu)^{-1} \in \sigma_{ei}(S_{\mu}),\eeq
and by Theorem ~\ref{Weyl}, \beq \sigma_{ei}(T_{\mu} - (T_{\mu} -
S_{\mu} )) \ = \
 \sigma_{ei}(T_{\mu}). \eeq

Applying Proposition  ~\ref{pee} shows that the result it true for
$i=5$ under the additional hypotheses.

\section{Further Notes and Remarks}

\vspace{0.3cm}
We note that Proposition ~\ref{closed} appeared for case  $ \sigma_{e1}$  in \cite{cross} (VII.2.3) and
that a similar but different generalisation of Weyl's theorem is proved in a lengthier argument through Theorems VII.2.15 and VII.2.3 of \cite{cross}.

\vspace{0.3cm}Other subsets of the spectrum of a linear operator have  also been
investigated for stability under perturbation, for example the
{\em Browder essential spectrum } defined by : \beq \sigma_b(T) \
:= \ \bigcup \{ \ \sigma(T+K)\ | \ \ TK=KT \ \ and \ \ K \ \ is \
\ compact \ \}. \eeq It is possible that  such investigations may
be extended to multivalued linear operators by the methods
employed in this work. More recently Sandovici, De Snoo and Winkler \cite{desnoo} have developed results for the ascent, descent, nullity and defect of linear relations.

\vspace{0.3cm} For simplicity, we have assumed  that the spaces on
which the relations are defined are complete, and that the
operators are closed. Fredholm properties are, however, stable
under more general conditions ( cf. Cross
\cite{cross} for the case $ \ \sigma_{e1} \ $ ). Thus, proofs for
$\ \sigma_{ei}, \ i = 1,2,3 \ $ do not necessarily require
assumptions of completeness. The index may not be stable under
perturbation,  though, and hence, generalisations which weaken
assumptions of completeness for $\ \sigma_{ei}, \ i = 4,5 \ $
would have to proceed with considerations similar to those applied
for the class of Atkinson relations  introduced in Wilcox \cite{wilcox} ( see also  L. Labuschagne \cite{lab} and V. M\"{u} ller-Horrig \cite{mh}).

\vspace{0.3cm}

%and Cross, Favini and Yakubov \cite{cfy} present apply spectral theory of linear relations to develope perturbation %theorems for a class of degenerate evolution equations.

\section{Acknowledgments}
The main body of this research was
undertaken while the author was based at the University of Cape
Town. The finalising of this manuscript was partially funded by
the National Research Foundation NPYY Programme for Young
Researchers [Grant Number 74223]. The author thanks R.W. Cross for the introduction to linear
relations in normed spaces and for catalysing this research by recommending references \cite{EdEv} and \cite{GuWe}.  The author also thanks M. M\"{o}ller for thoughtful questions on a preceding draft. Any errors herein are nevertheless my own.

\end{document}